\documentclass[11pt]{amsart}
\usepackage{graphicx}
\usepackage[ansinew]{inputenc}    % accents i altres

\newtheorem{theorem}{Theorem}
\theoremstyle{plain}

\newtheorem{definition}{Definition}
\newtheorem{example}{Example}

\newtheorem{lemma}{Lemma}

\newtheorem{proposition}{Proposition}
\newtheorem{remark}{Remark}

\numberwithin{equation}{section}

\begin{document}
\title{A characterization of quintic helices}
\author{J.V. Beltran, J. Monterde}
\address{Dep. de Geometria i
Topologia, Universitat de Val\`encia, Avd. Vicent Andr\'es
Estell\'es, 1, E-46100-Burjassot (Val\`encia), Spain}
\email{beltranv{@}uv.es, monterde{@}uv.es}

%\thanks{This work was partially supported by a Spanish MCyT grant
%BFM2002-00770.}
\date{February 16, 2005}

%\subjclass[2000]{Primary 53A04; Secondary 53C40, 53A05}


\begin{abstract}
A polynomial curve of degree $5$, $\alpha,$ is a helix if and only
if both $||\alpha'||$ and $||\alpha'\wedge \alpha''||$ are
polynomial functions.
\end{abstract}

\maketitle
\section{Introduction}

In \cite{FHMS} the authors study the notion of polynomial curves
which made a constant angle with a fixed direction. Curves with
constant slope are called there helices and we will use here the
same terminology in spite of the fact that in other contexts they
are also named generalized helices. The term ``helices'' is
reserved to curves in a cylinder with the same property.

We refer to the introduction of the cited paper and some other
papers like \cite{ChLM} or \cite{F} for the relationship between
such curves and some problems in the realm of computer-aided
design of curves and surfaces. The only fact we want to recall
here is that, for real applications, it seems clear that the
suitable curves are the quintic helices.

As it is said in \cite{FHMS}, any polynomial helix, $\alpha$, must
be a Pythagorean hodograph (PH) curve, i.e., $||\alpha'||^2$ is a
perfect square of a polynomial. Moreover, this condition is
sufficient in the cubical case: all PH cubics are helices.

As it is also said in the cited paper, (along the lines after
formula (11)) another necessary condition is the fact that
$||\alpha'\wedge \alpha''||^2$ (denoted there by $\rho^2$) must
also be a perfect square of a polynomial. The easiest way to see
this is the following: The argument that shows that a polynomial
helix must be PH, there applied to the tangent vector, ${\bf
\overrightarrow t} = \frac {\alpha'}{||\alpha'||}$, can also be
applied to the binormal vector, ${\bf \overrightarrow b} = \frac
{\alpha'\wedge\alpha''}{||\alpha'\wedge\alpha''||}$. (See the
first part of the proof of Th. 1 for details.) The consequence now
is that $||\alpha'\wedge\alpha''||$ is polynomial.

In this short paper we will go a little bit forward and show that
both conditions are sufficient in the quintic case. Moreover, we
will show an example of polynomial curve of degree $7$ verifying
both conditions but being not a helix.

\section{Spatial Pythagorean
hodograph curves}

We will use the quaternion representation of spatial PH curves.
Given a quaternion polynomial
$${\mathcal A}(t)= u(t)+{\bf i} v(t)+{\bf j} p(t)+{\bf k} q(t),$$
the product
$$\alpha'(t) ={\mathcal A}(t){\bf i}{\mathcal A}^\ast(t)$$
defines a spatial Pythagorean hodograph, $\alpha'$, whose
components are
\begin{eqnarray}\label{defcomponents}
x' &=& u^2+v^2-p^2-q^2,\nonumber\\
y' &=& 2(uq+vp),\\
z' &=& 2(vq-up),\nonumber
\end{eqnarray}
and such that $||\alpha'||^2=(x')^2+(y')^2+(z')^2
=(u^2+v^2+p^2+q^2)^2$.

In terms of the Hopf map (see \cite{ChLM}, theorem 4.2)
$$H:{\mathbb C}^2 \to {\mathbb R}^3$$ defined
by $H(z_1,z_2) = (|z_1|^2-|z_2|^2, 2 z_1 \bar{z}_2)$, and taking
\begin{equation}\label{z1z2} z_1(t) = u(t)+{\bf i}v(t),\qquad
z_2(t) = q(t)+{\bf i}p(t),
\end{equation}
the derivative of the curve can be written as
$$\alpha'(t) = H(z_1(t),z_2(t)).$$

\begin{definition} A curve $\alpha$ is said a Pythagorean hodograph curve of second class ($2$-PH curve)
if both $||\alpha'||$ and $||\alpha'\wedge \alpha''||$ are
polynomial functions.
\end{definition}

It is easy to check that 2-PH curves are examples of curves with a
rational Frenet-Serret frame (see \cite{WR}).

\begin{lemma}
The Frenet-Serret frame of any $2$-PH curve is made of rational
vectorial functions.
\end{lemma}
\begin{proof}
Simply recall that
$${\bf \overrightarrow t} = \frac {\alpha'}{||\alpha'||}, \quad
{\bf \overrightarrow b} = \frac
{\alpha'\wedge\alpha''}{||\alpha'\wedge\alpha''||}, \quad {\rm and
} \quad {\bf \overrightarrow n} =  {\bf \overrightarrow b}\wedge
{\bf \overrightarrow t}.$$
\end{proof}

\section{Characterization of $2$-PH curves}

Let $\alpha$ be a spatial PH curve whose tangent vector is defined
by the functions $u,v,p,q$ as in (\ref{defcomponents}) and let
$z_1, z_2$ be the associated complex functions as in (\ref{z1z2}).

\begin{proposition}\label{propo1}
$||\alpha'\wedge\alpha''||$ is a polynomial function if and only
if there is a complex polynomial function $z(t)$ and a real
polynomial function $\omega(t)$ such that
\begin{equation}\label{derivz1}
z_2^2 \left(\frac{z_1}{z_2}\right)' = \omega z^2.
\end{equation}
\end{proposition}

\begin{proof}
A straightforward computation shows that
$$||\alpha'\wedge\alpha''||^2 = 4 ||\alpha'||^2 ((u'q-uq'-v'p+vp')^2+(u'p-up'+
v'q-vq')^2).$$

Therefore, $||\alpha'\wedge\alpha''||$ is a polynomial function if
and only if $(u'q-uq'-v'p+vp')^2+(u'p-up'+ v'q-vq')^2$ is a
perfect square of a polynomial. Since both terms,
$u'q-uq'-v'p+vp'$ and $u'p-up'+ v'q-vq'$, are polynomial, we can
apply the well known result about Pythagorean curves, see \cite{F}
section 17.2: there is a polynomial function, $\omega(t)$, and a
complex polynomial function, $z(t)$, such that
$$(u'q-uq'-v'p+vp')+ {\bf i}(u'p-up'+ v'q-vq') = \omega  z^2.$$

An algebraic manipulation using the functions $z_1$ and $z_2$
defined in (\ref{z1z2}) allows to write the left hand member as
$$(u'+{\bf i} v')(q+{\bf i}p) - (u+{\bf i} v)(q'+{\bf i}p') =
z'_1 z_2 - z_1 z'_2 = z_2^2 \left(\frac{z_1}{z_2}\right)',$$ and
the statement follows.
\end{proof}

\begin{example}{\rm
Let us check this result in the two examples shown in \cite{FHMS}.

The first example is defined by the four quadratic polynomials
$$u(t) = t^2-3t,\quad v(t) = t^2-5t+10,\quad p(t) =
-2t^2+3t+5,\quad q(t) = t^2-9t+10.$$

Therefore, the complex functions $z_1(t) = (t^2-3t)+{\bf i}(
t^2-5t+10)$ and $z_2(t) = (t^2-9t+10)+{\bf i}(-2t^2+3t+5)$, verify
expression (\ref{derivz1}) for
$$
\omega(t) = 1,\qquad z(t) =\sqrt{1-7{\bf i}}\ (t-(1+2{\bf i})).$$
\medskip

The second example is defined by
$$\begin{array}{rcl}
u(t) &=& -19t^2+12 t+5,\qquad v(t) = -22 t^2+18 t+1,\\[2mm]
 p(t) &=&
15 t^2-12 t-1,\quad q(t) = -31 t^2+24t+3. \end{array}
$$

Now, expression (\ref{derivz1}) holds for
$$\omega(t) = 26(3 - 7 t + 3t^2),\qquad z(t) =\sqrt{-1+{\bf i}}.$$}
\end{example}

\section{$2$-PH curves of degree $5$.}

In this case $z_1(t)$ and $z_2(t)$  are quadratic polynomials and
the term on the left of the expression
$$z'_1 z_2 - z_1 z'_2 =\omega  z^2$$
is a polynomial of degree $2$. Therefore, there are two
possibilities for the pair $\omega(t), z(t)$. The first is
$\omega(t)$ be a quadratic function, and $z(t)$ a constant
function. The second, $\omega(t)$ be a constant
 and $z(t)$ a linear polynomial. As we will see, each possibility
correspond to one of the two classes of quintic helices studied in
\cite{FHMS}: the general helices and the monotone helices.

\medskip

We study first the case when $\omega(t)$ is constant, and without
loss of generality we can suppose that $\omega(t)=1$.

\begin{lemma}\label{lema1}
Monotone helices are characterized by a constant $\omega(t)$.
\end{lemma}
\begin{proof} The complex polynomials $z_1(t)$ and $z_2(t)$
are of degree less or equal than two. The first possibility is
that polynomials $z_1(t)$ and $z_2(t)$ are given by
$$\begin{array}{rcl}
z_1(t) &=& a(t-r_1)(t-r_2),\\
z_2(t) &=& b(t-r_3)(t-r_4),
\end{array}
$$
where $a,b,r_i\in{\mathbb C}$. An easy computation gives us that
$$\begin{array}{rcl} z'_1 z_2 - z_1 z'_2 &=&  ab
\left((r_1+r_2-r_3-r_4)\  t^2 + 2(r_3r_4-r_1r_2) t \right.\\
&& \left.+ (r_3+r_4)r_1r_2- (r_1+r_2)r_3r_4 \right).
\end{array}$$

If $\omega(t) = 1$ this expression is the square of a complex
polynomial function of degree $1$, $z(t) = m t + n$, if and only
if
$$\begin{array}{rcl}
ab(r_1+r_2-r_3-r_4) &=& m^2,\\
ab(r_3r_4-r_1r_2)&=& m n,\\
 ab((r_3+r_4)r_1r_2- (r_1+r_2)r_3r_4)&=& n^2.
\end{array}$$

From the first two equations we get
$$m = \pm\sqrt{ab}\sqrt{r_1+r_2-r_3-r_4}, \qquad n =
\frac{ab(r_3r_4-r_1r_2)}m.$$ Substituting in the last equation
$$r_1r_2r_3+r_1r_2r_4-r_1r_3r_4-r_2r_3r_4 = \frac
{r_1^2r_2^2-2 r_1r_2r_3r_4+r_3^2r_4^2}{r_1+r_2-r_3-r_4}.$$ After
some algebraic manipulation we can rewrite this equation as
$$(r_1-r_3)(r_2-r_3)(r_1-r_4)(r_2-r_4) = 0.$$

Therefore, $z'_1 z_2 - z_1 z'_2 = z^2$ if and only if $z_1(t)$ and
$z_2(t)$ share a linear factor. In this case
$gcd(z_1(t),z_2(t))\ne {\rm constant}$ which is the
characterization of monotone helices (see \cite{FHMS}, section
3.1). Indeed,
$${\rm gcd}(x',y',z') = |{\rm gcd}(u+{\bf i}v,p-{\bf i}q)|^2=
|{\rm gcd}(z_1,-{\bf i}z_2)|^2= |{\rm gcd}(z_1,z_2)|^2.$$

\medskip
The second possibility is that polynomials $z_1(t)$ and $z_2(t)$
are given by
$$\begin{array}{rcl}
z_1(t) &=& a(t-r_1)(t-r_2),\\
z_2(t) &=& b(t-r_3),
\end{array}
$$
where $a,b,r_i\in{\mathbb C}$. A similar analysis shows that $r_3
= r_1$ or $r_3= r_2$, and the same conclusion holds.

The last possibility is that polynomials $z_1(t)$ and $z_2(t)$ are
given by
$$\begin{array}{rcl}
z_1(t) &=& a(t-r_1)(t-r_2),\\
z_2(t) &=& b,
\end{array}
$$
where $a,b,r_i\in{\mathbb C}$. It is easy to check that this case
deals to a contradiction.
\end{proof}

\section{Characterization of quintic helices}

Let us recall first a description of PH quintics based on
quaternions, see \cite{FHMS}. A spatial quintic helix is defined
by a quadratic polynomial
$${\mathcal A}(t) ={\mathcal A}_0 + {\mathcal A}_1 t + {\mathcal
A}_2 t^2,$$ with quaternion coefficients
$$\begin{array}{rcl}
{\mathcal A}_0 &=& a + a_x{\bf i}+ a_y{\bf j}+ a_z{\bf k},\\[2mm]
{\mathcal A}_1 &=& b + b_x{\bf i}+ b_y{\bf j}+ b_z{\bf k},\\[2mm]
{\mathcal A}_2 &=& c + c_x{\bf i}+ c_y{\bf j}+ c_z{\bf k}.
\end{array}
$$
In terms of the functions $u,v,p,q$:  $u(t)=a+b t+c t^2$,
$v(t)=a_x+b_x t+c_x t^2$, $p(t)=a_y+b_y t+ c_y t^2$ and
$q(t)=a_z+b_z t+c_z t^2$.

\begin{lemma}\label{lema2}
Let  $z_1(t) = u(t)+{\bf i}v(t)$ and  $z_2(t) = q(t)+{\bf i}p(t)$
be the quadratic polynomials of a quintic PH curve defined by
three quaternions $\{{\mathcal A}_i\}$, $i=0,1,2$. If $z(t)$ in
(\ref{derivz1}) is constant  then ${\mathcal A}_1 = c_0 {\mathcal
A}_1 +c_2 {\mathcal A}_2 $, for suitable real scalars $c_0,c_2$.
\end{lemma}

\begin{remark}{\rm
In \cite{FHMS} the authors use a Bézier quadratic polynomial
$${\mathcal A}(t) ={\mathcal A}_0 (1-t) ^2 + 2{\mathcal A}_1 t(1-t) + {\mathcal
A}_2 t^2.$$ We have used here the usual basis of polynomials
instead of the Bernstein basis because computations are easier.
The statement of the previous Lemma remains true for Bézier
quaternion coefficients due to just a change of basis.}
\end{remark}

\begin{proof} Let us suppose that
$$\omega(t) = m_0 + m_1 t + m_2 t^2,\quad  z(t) = e^{{\bf i}\theta},$$
where, $m_0, m_1, m_2, \theta\in{\mathbb R}$ and without loss of
generality, we assume that  $|z| = 1$.

We use now Proposition \ref{propo1} and compute the expression
$$z'_1 z_2 - z_1 z'_2 =\omega z^2.
$$
The real part of the left hand term can be written as
$$
\begin{array}{rl}
&u'q-uq'- v'p+vp' = (a_x b_y - a b_z + a_z b - a_y  b_x)\\
 &\qquad + 2 ( a_x c_y - a c_z - a_y  c_x + a_z c) t + (b_z c - b_y c_x + b_x c_y - b c_z)  t^2
, \end{array}
$$
and the imaginary part as
$$
\begin{array}{rl}
&u'p-up'+ v'q-vq' = (a_y b + a_z b_x - a
b_y - a_x b_z) \\[2mm]
& \qquad + 2(a_y c + a_z c_x - a c_y - a_x c_z) t + (b_y c + b_z
c_x - b c_y - b_x c_z)t^2. \end{array}
$$

Analogously, the real part of the right hand term can be written
as
$$(m_0 + m_1 t + m_2 t^2) \cos(2\theta),$$
and the imaginary part as
$$(m_0 + m_1 t + m_2 t^2) \sin(2\theta).$$

Therefore, the condition $z'_1 z_2 - z_1 z'_2 =\omega z^2$ can be
translated into a set of six equations. By equating the
coefficients of $t$ we can deduce that
$$\theta = \frac 12\arctan(\frac{a_y c + a_z c_x - a c_y - a_x c_z}{a_z c - a_y c_x + a_x c_y -
a c_z})$$ and $$m_1 = 2\sqrt{(a_y c + a_z c_x - a c_y - a_x c_z)^2
+ (a_z c - a_y c_x + a_x c_y - a c_z)^2}.$$

Substituting these values into the other four equations and
solving the resulting linear system we obtain
$${\mathcal A}_1 = \frac{2m_0}{m_1} {\mathcal A}_0 + \frac{2m_2}{m_1} {\mathcal
A}_2.$$
\end{proof}

\begin{theorem}
A quintic polynomial curve is a  helix if and only if it is a
$2$-PH curve.
\end{theorem}

\begin{proof}
(Necessary conditions)

Let us recall that if a curve is a helix, then, not only the
tangent vector  makes a constant angle with the axis, but also the
binormal vector, see the classical references \cite{dC},\cite{St}.
In fact, if $\overrightarrow{u}$ is a unitary vector that
determines the axis of the helix then
$$<\overrightarrow{u},\overrightarrow{t}> = c, \qquad
<\overrightarrow{u},\overrightarrow{b}> = \sqrt{1-c^2},$$ where
$c\in{\mathbb R}$ is a constant, $\overrightarrow{t}=\frac
{\alpha'}{||\alpha'||}$ is the tangent vector and
$\overrightarrow{b}=\frac {\alpha'\wedge \alpha''}{||\alpha'\wedge
\alpha''||}$ is the binormal vector of the curve.

The previous expressions are equivalent to
$$<\overrightarrow{u},\alpha'> = c \ ||\alpha'||,\qquad
<\overrightarrow{u},\alpha'\wedge \alpha''> = \sqrt{1-c^2}\
||\alpha'\wedge \alpha''||.$$

If our curve $\alpha$ is polynomial then it is a Pythagorean
hodograph and also the norm of $||\alpha'\wedge \alpha''||$ is
polynomial, indeed
$$||\alpha'||=\frac{1}{c}\ <\overrightarrow{u},\alpha'>  \ ,\qquad
||\alpha'\wedge \alpha''||=\frac{1}{\sqrt{1-c^2}}\
<\overrightarrow{u},\alpha'\wedge \alpha''>  \ .$$

\bigskip

(Sufficient conditions) If $||\alpha'\wedge \alpha''||$ is
polynomial then by Proposition \ref{propo1} we know that $z_2^2
\left(\frac{z_1}{z_2}\right)' = \omega z^2$. In the quintic case,
the only possibilities are $w$ constant or $z$ constant. If
$\omega$ is constant then by Lemma \ref{lema1} the curve is a
monotone helix.

If $z$ is constant, then by Lemma \ref{lema2} the quaternions
defining the curve are linear dependents and by Proposition $1$ in
\cite{FHMS} we know that the curve is a helix.
\end{proof}

\begin{remark}{\rm
In higher dimensions it is possible to find $2$-PH curves being
not helices. For example:
$$\alpha(t) = (-3t + t^3 + \frac{t^5}{5} + \frac{t^7}{21}, 3 t^2 - \frac{t^4}2, -2 t^3)$$
is a polynomial curve of degree $7$ verifying
$$||\alpha'|| = \frac 13( 9 + 9 t^2 + 3 t^4 + t^6),\quad
||\alpha'\wedge\alpha''|| = 2(1 + t^2)( 9 + 9 t^2 + 3 t^4 + t^6)$$
but $$\frac{\tau}{\kappa} = \frac{-9+ 9 t^4 + 2 t^6}{9(1 +
t^2)^2},$$ so, it does not satisfy the Lancret's theorem (see for
example \cite{St}) and the curve is not a generalized helix.}
\end{remark}

\end{document}